\documentclass[11pt]{amsart}

\usepackage{amscd}
\usepackage{amssymb}
\usepackage{boxedminipage}
\usepackage{enumerate}

\theoremstyle{plain} \newtheorem{theorem}{Theorem}[section]
\theoremstyle{plain} 
\theoremstyle{plain} 
\theoremstyle{plain}

\newcommand{\nr}{\refstepcounter{theorem}  
                   \noindent {\thetheorem .}}

\newcommand{\eks}{\medskip \noindent {\it Example \nr} }
\newcommand{\eksfin}{\medskip}
\newcommand{\rem}{\medskip \noindent {\it Remark \nr} }

\newcommand{\remfin}{\medskip}

\newcommand{\al}{\alpha}
\newcommand{\ep}{\epsilon}

\newcommand{\proj}[1]{{\bf P}(#1)}

\newcommand{\dpu}{v}
\newcommand{\pdu}{v}

\newcommand{\modv}[1]{{#1}\text{-{mod}}}

\newcommand{\ome}{\omega_E}
\newcommand{\coh}{\text{coh}}

\newcommand{\op}{{\mathcal O}}
\newcommand{\opu}{{\mathcal O}_{\proj{U}}}

\newcommand{\gF}{{\mathcal F}}

\newcommand{\gE}{{\mathcal E}}

\newcommand{\gR}{{\mathcal R}}

\newcommand{\gL}{{\mathcal L}}

\newcommand{\cfr}{{-cF}}
\newcommand{\ecofree}{E\!\cfr}

\newcommand{\Hom}{\text{Hom}}

\newcommand{\im}{\text{im}\,}

\newcommand{\ind}{\text{ind}}
\newcommand{\diag}{\text{diag}}

\newcommand{\sus}{\subseteq}
\newcommand{\larr}{\longrightarrow}
\newcommand{\pil}{\rightarrow}
\newcommand{\lpil}{\larr}
\newcommand{\inpil}{\hookrightarrow}

\newcommand{\mto}[1]{\stackrel{#1}\longrightarrow}

\newcommand{\iso}{\cong}
\newcommand{\te}{\otimes}

\newcommand{\lvi}{\langle}
\newcommand{\hvi}{\rangle}

\newcommand{\pn}{{\bf P}^r}

\newcommand{\go}{{\mathcal O}}

\newcommand{\hele}{{\bf Z}}

\begin{document}
\title {Exterior algebra resolutions arising from homogeneous bundles}
\author { Gunnar Fl{\o}ystad}
\address{ Matematisk Institutt\\
          Johs. Brunsgt. 12 \\
          5008 Bergen \\
          Norway}   
        
\email{ gunnar@mi.uib.no }

\maketitle

\section*{Introduction.}

When $R$ is a commutative ring, the minimal free resolution of a
map $R^a \pil R^b$ and symmetric and skew-symmetric maps
$R^a \pil R^a$, under suitable generality conditions,
are well known and have been developed by a series of
authors. See \cite[A2.6]{Ei} for an overview.

In this note we do an analog for the exterior algebra
$E = \oplus \wedge^i V$ on a finite dimensional vector space $V$ 
and general graded maps $E^a \pil E(1)^b$, and
general graded symmetric and skew-symmetric maps $E^a \pil E(1)^a$. 
Since $E$ is 
both a projective and injective $E$-module, by taking a free (projective)
and cofree (injective) resolution of such maps, there is associated
an unbounded acyclic complex of free $E$-modules, called a Tate resolution,
see \cite{EFS} or \cite{Fl}.
Via the Bernstein-Gel'fand-Gel'fand (BGG) correspondence, this corresponds
to a complex of coherent sheaves on the projective space $\proj{V^*}$. 
We show that in all the cases above (with the dimension of $V$ not too low),
this complex actually reduces to a
coherent sheaf. We describe these coherent sheaves and also describe
completely the Tate resolutions.

These descriptions turn out to be simpler to work
out than in the corresponding commutative case, and, maybe at first surprising,
the descriptions are also more geometric.

In fact not only are we able to describe the Tate resolutions
and coherent sheaves
associated to the maps stated above, but, using the theory
of representations of reductive groups, we are able to describe the Tate
resolutions and coherent sheaves associated to vast larger classes of 
natural maps $E^a \pil E(r)^b,$ something which would have required
considerably more effort for commutative rings.

There is only one catch related to all our descriptions. We must assume that
the dimension of $V$ is not too small compared to $a$ and $b$. For instance
for a general map $E^a \pil E(1)^b$ we must assume that the dimension of $V$
is  $\geq a+b-1$.
In case the dimension of $V$ is smaller than this, the nature of the problem changes, and we do not investigate this case.



\section{ Tate resolutions and projections.}

Let $V$ be a finite dimensional vector space over a field $k$. 
Put $W = V^*$ and
set $\dpu = \dim \proj{W}$. Let $E(V)$ be the graded exterior algebra 
$\oplus_{i=0}^{\dpu + 1}\wedge^i V$
where $V$ has degree $-1$,
(because we consider $W$ to have degree 1). Let $\ome = \Hom_k(E,k)$
be the graded dual, which we consider as a left $E$-module.
(As such $\ome \iso E(-\dpu-1)$.) It is the injective hull of $k$.

\subsection{Tate resolutions.}
A {\it Tate resolution} is an (unbounded) acyclic complex $T$ with components
\[ T^p = \oplus_{i \in \hele} \ome(-i) \te_k V^p_i \]
of finite rank.
Note that a Tate resolution is completely determined, up to homotopy, by each
differential $d^p$ since $T^{\leq p}$ is a projective resolution of 
$\im d^p$ and $T^{> p}$ is an injective resolution of $\im d^p$.

By \cite{EFS} or \cite{Fl},
to each coherent sheaf $\gF$ on $\proj{W}$ there is associated a Tate
resolution $T(\gF)$ whose terms are
\[ T(\gF)^p = \oplus_{i = 0}^\dpu \ome(i-p) \te_k H^i \gF(p-i). \]
(In particular we see that $T(\gF)^p$ is  $\ome(-p) \te_k H^0 \gF(p)$ for
$p \gg 0$.) The maps 
\[ \ome(i-p) \te_k H^i \gF(p-i) \lpil \ome(i-p-1) \te_k H^i \gF(p+1-i) \]
are determined by the maps in degree $p+1-i$ which is the natural map 
\[ W \te_k H^i \gF(p-i) \lpil H^i \gF(p+1-i). \]

Conversely, given a left module $N = \oplus_{i \in \hele} N_i$
over $E$ we get associated a {\it complex} of coherent sheaves 
\[ G(N)^{\sim} :\,\, \cdots \pil \opu(i) \te_k N_i \pil \opu(i+1) \te_k N_{i+1} \pil \cdots . \]
In this way we can to each Tate resolution $T$ associate a complex of coherent
sheaves on $\proj{W}$ by using this construction on $\ker d^p$.
(All the $G(\ker d^p)^{\sim}$ become isomorphic in the derived
category of coherent sheaves on $\proj{W}$. If $T = T(\gF)$
then $G(\ker d^p)^{\sim}$ is a complex isomorphic to $\gF$ in this derived
category.)

\rem If $K^\circ(\ecofree)$ is the homotopy category of Tate resolutions,
and $D^b(\coh/\proj{W})$ is the derived category of coherent sheaves
on $\proj{W}$,
the Bernstein-Gel'fand-Gel'fand correspondence says that
there is an equivalence of categories
\[ K^\circ(\ecofree) \iso D^b (\coh/\proj{W}) . \]
See \cite{Fl} or originally \cite{BGG}.
\remfin

\subsection{Projections}
Let $U \sus W$ be a subspace and let $E(U^*) = \oplus \wedge^i U^*$. The center
of the projection $p: \proj{W} \dashrightarrow \proj{U}$ is the linear
subspace $\proj{W/U}$. Suppose Supp$\,\gF \cap \proj{W/U} = \emptyset$. 
According to \cite{Fl},
\[ \Hom_{E(V)}(E(U^*), T(\gF)) \]
is also a Tate resolution (for $E(U^*)$) and it is the Tate resolution
associated to $p_* \gF$. 
Note that this gives
\[ T(p_*\gF)^p = \oplus_{i = 0}^{\dim \proj{U}}
\omega_{E(U^*)}(i-p) \te_k H^i \gF(p-i). \]
\medskip

If $Y \pil W$ is a linear map, then let $U$ be the image in $W$.
Thus we get a projection and and embedding
\[ \proj{W} \dashrightarrow \proj{U} \mto{i} \proj{Y}. \]
In this case 
\[ \Hom_{E(V)}(E(Y^*), T(\gF)) \]
is the Tate resolution of $i_* p_* \gF$.

\section{Tate resolutions arising from general maps.}

Let $A$ and $B$ be finite dimensional vector spaces. A general map
\[ A^* \lpil V \te B \]
comes from the generic map
\[A^* \lpil (A^* \te B^*) \te B \]
composed with a general map $A^* \te B^* \pil V$. It
gives us a general morphism
\[ \ome \te A^* \lpil \ome(-1) \te B \]
which again gives a Tate resolution by taking projective and injective
resolutions. This again is associated to a complex of coherent sheaves
which we now show is a coherent sheaf. Let $a$ and $b$ be the dimensions
of $\proj{A}$ and $\proj{B}$.

\begin{theorem} Given a general morphism
\begin{equation} \ome(a) \te A^* \mto{d} \ome(a-1) \te B \label{2genavb}
\end{equation}
coming from a general surjection $A^* \te B^* \pil V$ and suppose 
$\dpu \geq a + b$.
Let $T$ be the associated Tate resolution with (\ref{2genavb}) in components 
$0$ and $1$.

a. $T$ is the Tate resolution associated to $p_* \gL$ where 
$\gL$ is the line bundle $\op_{\proj{A} \times \proj{B}}(-2,a) 
\te \wedge^{a+1} A$
on the Segre embedding of $\proj{A} \times \proj{B}$ in $\proj{A \te B}$
and $p : \proj{A \te B} \dashrightarrow \proj{W}$ is the projection.

b. The Tate resolution has terms
\[ T^p = \begin{cases} \ome(-p) \te S^{p-2}A \te S^{p+a}B \te \wedge^{a+1} A &
,\,\, p \geq 2, \\
\ome(a+b+r) \te S^{b+r-1} A^* \te S^{r-1} B^* \te \wedge^{b+1} B^* & 
,\,\, p = -r < 0
\end{cases} \]
\end{theorem}

\begin{proof}
We find
\begin{eqnarray*}  
H^a \gL(-a) & = & \wedge^{a+1}A \te 
H^a \go_{\proj{A} \times \proj{B}}(-2-a,0)\\
   & = & \wedge^{a+1}A \te A^* \te \wedge^{a+1}A^* = A^*
\end{eqnarray*}

\begin{eqnarray*}
H^a \gL(-a+1) & = & \wedge^{a+1}A \te 
                    H^a \go_{\proj{A} \times \proj{B}} (-1-a,1)\\
            & = & \wedge^{a+1}A \te \wedge^{a+1}A^* \te B = B
\end{eqnarray*}
and also $H^i \gL(-i)$ and $H^i \gL(-i+1)$ are zero for $i \neq a$.
The map
\[ (A \te B) \te H^a \gL(-a) \pil H^a \gL(-a+1) \]
is clearly the generic one and hence we get part a. 
Part b. follows of course by considering the cohomology of
$\gL$.
\end{proof}

\rem Consider the generic morphism 
\[ \op_{\proj{A \te B}} \te A^*  \pil \op_{\proj{A \te B}}(1) \te B. \]
The rank strata of this morphism are the orbits
in $\proj{A \te B}$ of $GL(A) \times GL(B)$. The stratum of lowest rank
one, corresponds to the closed orbit of $GL(A) \times GL(B)$ which is 
the Segre embedding $\proj{A} \times \proj{B} \inpil \proj{A \te B}$.
\remfin

\section{The Tate resolution of a general symmetric map.}

Given now a symmetric map
\[ A^* \lpil V \te A. \]
It comes from the natural map
\[ A^* \pil S^2 A^* \te A \]
composed with a map $S^2 A^* \pil V$.
We then get a symmetric morphism
\[ \ome \te A^* \pil \ome(-1) \te A. \]

\begin{theorem}
Given a general symmetric morphism
\begin{equation} \ome(a) \te A^* \mto{d} \ome(a-1) \te A 
\label{3symmavb} \end{equation}
coming from a general surjection $S^2 A^* \pil V$ with $\dpu \geq 2a$.
Let $T$ be the associated Tate resolution with (\ref{3symmavb}) as the
components in degree $0$ and $1$.
a. $T$ is the Tate resolution associated to $p_* \gL$ where
$\gL$ is the line bundle 
$\op_{\proj{A} \times \proj{A}}(-2,a) \te \wedge^{a+1} A$
on the Segre embedding of $\proj{A} \times \proj{A}$ in $\proj{A \te A}$
and $p : \proj{A \te A} \dashrightarrow \proj{S^2 A} \dashrightarrow 
\proj{W}$ is the projection.

b. The Tate resolution has terms
\[ T^p = \begin{cases} \ome(-p) \te S^{p-2}A \te S^{p+a}A \te \wedge^{a+1} A &
,\,\, p \geq 2, \\
\ome(2a+r) \te S^{a+r-1} A^* \te S^{r-1} A^* \te \wedge^{a+1} A^* & 
,\,\, p = -r < 0
\end{cases}\]
\end{theorem}

\begin{proof}

This follows from the previous theorem by verifying that the
Segre embedding $\proj{A} \times \proj{A} \inpil \proj{A \te A}$
does not intersect the center of the projection 
$\proj{A \te A} \dashrightarrow \proj{S^2 A}$.


\end{proof}

We now wish to describe the Tate resolution and coherent sheaf associated
to a general skew-symmetric morphism. From the preceding one might thing 
that this is quite analogous to what we have done for symmetric maps.
However it is quite different due to the fact that the projection 
center of $\proj{A \te A} \dashrightarrow \proj{\wedge^2 A}$ 
intersects the Segre embedding $\proj{A} \times \proj{A} \inpil
\proj{A \te A}$.
We must therefore take another approach which will give us a much more 
powerful understanding of what we have just done.

\section {Homogeneous bundles on homogeneous spaces.}

We present here some facts about induced homogeneous bundles 
on homogeneous varieties. In particular we consider their cohomology.
We base ourselves on \cite{Ja}.

Let $G$ be a reductive algebraic group over a field $k$ of characteristic
zero, with Borel
group $B$ and maximal torus $T$.
We let $X(T)$ be the characters of $T$ and $R$ the root system of $G$.
Corresponding to the choice of $B$ there is a positive root system $R^+$
and let $S$ be the associated simple roots.


Given a parabolic subgroup $P$ of $G$, we let $\modv{P}$ be the category
of rational representations of $P$. For each $P$-module $M$ we get
a $G$-equivariant locally free sheaf $\gL(M)$ of rank $\dim_k M$
on the homogeneous variety $G/P$ \cite[I.5.8]{Ja}.
If $N$ is another $P$-module, then $\gL(M \te_k N) =
\gL(M) \te_{\go_{G/P}} \gL(N)$.

To such a parabolic $P$ is associated a subset $I \sus S$ and the 
characters $X(P)$ are all $\mu \in X(T)$ such that $\lvi \mu, \check{\al} 
\hvi = 0$ for all $\alpha \in I$. For such $\mu$ we get line bundles
$\gL(\mu)$ on $G/P$. By \cite[II.4.6]{Ja}, $\gL(\mu)$ is ample
(in fact very ample) on $G/P$ if
$\lvi \mu, \check{\al} \hvi > 0$ for all $\alpha \not \in I$. Thus
we get an embedding
\[ G/P \inpil \proj{H^0(G/P, \gL(\mu))} \]
with $\op_{G/P}(1) = \gL(\mu)$.

Let 
\[ \ind^G_P : \modv{P} \pil \modv{G} \]
be the induction functor \cite[I.3.3]{Ja} and let
$R^n \ind^G_P$ be its higher derived functors.
Then $H^n (G/P, \gL(M)) = R^n \ind^G_P M$
by \cite[I.5.10]{Ja}. If $\lambda \in X(T)$, to simplify notation, let 
$R^n \ind^G_B \lambda = H^n(\lambda)$. This is an irreducible $G$-module
whose character is determined by the Borel-Bott-Weil theorem 
\cite[II.5.5]{Ja} and the Weyl character formula \cite[II.5.10]{Ja}.

Associated to the inclusions $B \sus P \sus G$ there is by \cite[I.4.5]{Ja}
a spectral sequence
\begin{eqnarray} E_2^{n,m} = R^n \ind^G_P \circ R^m \ind^P_B N \Rightarrow 
R^{n+m} \ind^G_B N. \label{4spect}
\end{eqnarray}
If now $\lambda \in X(T)$ is such that $\lvi \lambda, \check{\al} \hvi 
\geq 0$ for all $\alpha$ in $I$, then
$R^n \ind^P_B \lambda = 0$ for $n > 0$ by Kempf vanishing \cite[II.4.5]{Ja}.
In this case we therefore get by (\ref{4spect}) that
\begin{eqnarray*}
H^n(\lambda) & = & R^n \ind^G_P (\ind^P_B \lambda) \\
             & = & H^n (G/P, \gL(\ind^P_B \lambda)).
\end{eqnarray*}
Considering the $G$-equivariant induced sheaf $\gL(\ind^P_B \lambda)$ on
$G/P \inpil \proj{H^0(\mu)}$ we are thus able to determine all the 
cohomology groups
\begin{eqnarray*}
H^n (G/P, \gL(\ind^P_B \lambda)(r)) & = & H^n(G/P, \gL(\ind^P_B \lambda) 
\te \gL(r \mu)) \\
& = & H^n (G/P, \gL(\ind^P_B (\lambda + r \mu))) = H^n (\lambda + r \mu).
\end{eqnarray*}

In particular, for the induced homogeneous bundle 
$\gL = \gL(\ind^P_B \lambda)$
on $\proj{H^0(\mu)}$, the terms of the Tate resolution $T$ of $\gL$
\[ T^p = \oplus_{i = 0}^{\dim G/P} \ome(i-p) \te H^i \gL(p-i) \]
may be determined.

\eks
Let $G = GL(W)$ and let $e_0, \ldots, e_\dpu$ be a basis for $W$.
If $T$ is the diagonal matrices we let $\ep_i$ be the character sending
$\diag(t_0, \ldots, t_\dpu)$ to $t_i$. We let $B$ be the lower
triangular matrices. The positive roots are then
$\ep_i - \ep_j$ where $i < j$ and the positive simple roots are
$\al_i = \ep_{i-1} - \ep_i$ for $i = 1 \ldots \dpu$.
Let $P$ be the parabolic subgroup corresponding to $\al_2, \ldots ,\al_\pdu$,
i.e. $P$ consists of all matrices
\[ \begin{pmatrix} * & 0 & \cdots & 0 \\
                   * & * & \cdots & * \\
                   \vdots & & & \\
                   * & * & \cdots & * 
\end{pmatrix}. \]
Letting $U = (e_1, \ldots, e_\pdu)$ there is an exact sequence of 
$P$-modules 
\[ 0 \pil U \pil W \pil (e_0) \pil 0. \]
Here $U$ is the irreducible representation of $P$ with highest weight 
$\ep_1$ and thus $U = \ind^P_B \ep_1$.
On $G/P = \proj{W}$ there is an exact sequence
\[ 0 \pil \gL(U) \pil \gL(W) \pil \gL(\ep_0) \pil 0. \]
Since the Picard group of $G/P$ is generated by $\gL(\ep_0)$ we have
$\gL(\ep_0) = \opu(1)$. Since $W$ is a $G$-module, 
$\gL(W) = \opu \te_k W$. Thus $\gL(U) = \Omega_{\proj{W}}(1)$.

Given now a partition 
\[ {\bf{i}} : i_1 \geq i_2 \geq \cdots \geq i_\pdu \]
where $i_j$ are integers.  The Schur bundle 
$S_{\bf{i}}(\Omega_{\proj{W}}(1))$ is then the induced bundle
$\gL(S_{\bf{i}} U)$ where $S_{\bf{i}}$ is the Schur functor, \cite[A2.5]{Ei}. 
Since $S_{\bf{i}} U$ is an irreducible representation
with highest weight $i_1 \ep_1 + \cdots + i_\pdu \ep_\pdu$ we get
$S_{\bf{i}} U = \ind^P_B ( \sum_{j=1}^n i_j \ep_j)$.
It then follows that
\begin{eqnarray*} H^r S_{\bf{i}}(\Omega_{\proj{W}}(1))(p-r) & = & 
H^r (\proj{W}, \gL(S_{\bf{i}} U) \te \gL ((p-r) \ep_0)) \\
 & = & H^r (G/P, \gL(\ind^P_B ((p-r) \ep_0 + \sum_{j = 1}^n i_j \ep_j))) \\
 & = & H^r ((p-r)\ep_0 + \sum_{j=1}^n i_j \ep_j).
\end{eqnarray*}

By the Borel-Bott-Weil theorem this can be calculated and the answer
turns out to be the following.

Let $h = h(p)$ be such that $i_h > p \geq i_{h+1}$ and let
${\bf{i}}(p)$ be the partition
\[ i_0 - 1 \geq \cdots \geq i_h - 1 \geq p \geq i_{h+1} \geq 
\cdots \geq i_\pdu. \]
It then follows that
\[ H^r ((p-r) \ep_0 + \sum_{j=1}^n i_j \ep_j) = \begin{cases} 
  H^0 (\sum_{i = 0}^{\dpu} {\bf{i}}(p)_j \ep_j) &,\,  r = h(p) \\
 0 &,\,  r \neq h(p)
\end{cases}. \]

In particular the Tate resolution $T$ of $S_{\bf{i}}(\Omega_{\proj{W}}(1))$
has "pure" terms
\[ T^p = \ome(h(p) - p) \te S_{{\bf{i}}(p)} W. \]

\rem More generally one can show that if $\mu = \omega_\al$ is a fundamental
weight corresponing to a short root and $P_{S-\{\al\}}$ is the corresponding
parabolic group, then all the induced bundles 
$\gL(\ind_B^{P_{S-\{\al\}}} \lambda)$ have "pure" Tate resolutions in
the above sense.
\remfin

\section{ Tate resolutions arising from general skew-symmetric maps.}

Suppose given a skew-symmetric map
\[ A^* \pil V \te A. \]
Such a map comes from the natural map
\[ A^* \pil \wedge^2 A^* \te A \]
and a map $\wedge^2 A^* \pil V$.
It gives rise to a skew-symmetric morphism
\[ \ome \te A^* \pil \ome(-1) \te A. \]

To describe the Tate resolution and associated coherent sheaf on 
$\proj{W}$ consider the Pl\"ucker embedding of the Grassmann of lines
in $\proj{A}$
\[ G(A,2) \inpil \proj{\wedge^2 A}.\]
Let 
\[ 0 \pil \gR \pil \op_G \te A \pil \gE \pil 0 \]
be the tautological sequence (where the rank of $\gE$ is two).

\begin{theorem}
Given a general skew-symmetric morphism
\begin{equation} \ome(a-1) \te A^* \pil \ome(a-2) \te A \label{5ssavb}
\end{equation}
arising from a general surjection $\wedge^2 A^* \pil V$ and suppose
$\dpu \geq 2(a-1)$. Let $T$ be the associated Tate resolution
with components in degree $0$ and $1$ given by (\ref{5ssavb}).

a. $T$ corresponds to the coherent sheaf $p_* \gL$ where
$\gL$ is the bundle $(S^a \gE)(-2) \te \wedge^{a+1}A^*$ on
$G(A,2) \inpil \proj{\wedge^2 A}$ and
$p : \proj{\wedge^2 A} \dashrightarrow \proj{W}$ is the projection.
 
b. The components of $T$ are given by
\[ T^p = \begin{cases} \ome(-p) \te S_{a+p-2, p-2} A \te \wedge^{a+1} A 
 &,\,\, p \geq 2 \\
  \ome(2a+r-2) \te S_{a+r-1,r-1} A^* \te \wedge^{a+1} A^* &,\,\, p = -r < 0
\end{cases}. \] 
\end{theorem}

\begin{proof} Choose a basis $e_0, \ldots, e_a$ for $A$ and let $B$ be
the Borel subgroup of lower triangular matrices in $GL(A)$. Let $(e_i)$ be the
weight space of $\ep_i$. Then $G(A,2) = GL(A)/P$ where $P$ is the
maximal parabolic subgroup associated to the roots $\al_i = 
\ep_{i-1} - \ep_i$, $i \neq 2$, i.e. $P$ consists of matrices 
\[ \begin{pmatrix} * & * & 0 &\cdots & 0 \\
                   * & * & 0 &\cdots & 0 \\
                   * & * & * & \cdots & * \\
                    \vdots & & & &\\
                   * & * & * & \cdots & * 
\end{pmatrix}. \]

Note that $\dim G(A,2) = 2(a-1)$ and that by the Pl\"ucker embedding
$G(A,2) \inpil \proj{\wedge^2 A}$ the canonical bundle
$\op_G(1)$ is $\gL(\ep_0 + \ep_1)$, where $\ep_0 + \ep_1$ is the highest
weight of $\wedge^2 A$.

There is a sequence of $P$-modules 
\[ 0 \pil (e_2, \ldots, e_a) \pil A \pil (e_0, e_1) \pil 0 \]
and the tautological sequence on $G(A,2)$ becomes
\[ 0 \pil \gL((e_2, \ldots, e_a)) \pil \op_G \te A 
\pil \gL((e_0, e_1)) \pil 0. \]
Since $(e_0, e_1) = \ind_B^P \ep_0$ as $P$-modules we get
\[ S^a \gE = S^a \gL((e_0,e_1)) = \gL(S^a(e_0,e_1)) = \gL(\ind_B^P(a \ep_0)) \]
and 
\[ \gL := (S^a \gE)(-2) \te \wedge^{a+1} A^* = 
\gL( \ind_B^P ( (a-1)\ep_0 - \ep_1 + \sum_{2}^a \ep_i )). \]
By the Borel-Bott-Weil theorem we find that
\[ H^{a-1} \gL(-a+1) = H^{a-1}( -a \ep_1 + \sum_2^a \ep_i) 
  = H^0 (-\ep_a) = A^* \]
and 
\[ H^{a-1} \gL(-a+2) = H^{a-1}(\ep_0 - (a-1)\ep_1 + \sum_2^a \ep_i)
                    = H^0 (\ep_0) = A. \]
Since we find $H^i \gL(-i) = 0$ and $H^i \gL(-i+1) = 0$ for $i \neq a-1$,
part a. follows.

We find that 
\begin{equation*}
H^0 \gL(p) = H^0 ((a+p-1) \ep_0 + (p-1)\ep_1 + \sum_2^a \ep_i)
           = S_{a+p-2, p-2} A \te \wedge^{a+1} A, \quad p \geq 2
\end{equation*}
and by Borel-Bott-Weil that
\begin{eqnarray*} H^{2a-2} \gL(2-2a-r) & = & 
\wedge^{a+1} A \te H^{2a-2} ((-a-r)\ep_0 + (-2a -r) \ep_1) \\
         & = & \wedge^{a+1} A \te (\wedge^{a+1} A^*)^{\te 2} \te 
         S_{a+r-1, r-1} A^* \\
         & = & S_{a+r-1, r-1} A^* \te \wedge^{a+1} A^*, \quad p = -r < 0
\end{eqnarray*}
and also that $H^i \gL(p-i) = 0$ for all other $i$ and $p$.
Hence b. follows.

\end{proof}

\section {Some further examples.}

\eks Consider the map 
\[ \ome(1) \te V \mto{\phi} \ome(-1) \te V^* \]
corresponding to the map
$ \wedge^2 V^* \te V \pil V^* . $
Here we do not get associated a coherent sheaf to the Tate resolution
$T(\phi)$ but instead a complex, the truncated Koszul complex
\begin{equation} \op_{\proj{W}}(-\pdu) \te \wedge^\pdu W \pil \ldots 
\pil \op_{\proj{W}}(-2) \te \wedge^2 W \pil \op_{\proj{W}}(-1) \te W 
\label{trunkos6} \end{equation}

In fact, let $M$ be $\oplus_{i=1}^{\pdu} \wedge^i W$ which is a subquotient of 
$\ome$. Then we see that (\ref{trunkos6}) is $G(M)^{\sim}$.

Since $\ome(1) \te V$ is the projective cover of $M$ and 
$\ome(-1) \te W$ is the injective hull, we get $M = \im \phi$ and
so $T(\phi)$ corresponds to (\ref{trunkos6}).
\eksfin

\eks Suppose $\dpu$ is even. Let $\wedge^2 W \pil k$
be a non-degenerate symplectic form.
Then the map $\wedge^{\dpu} W \pil k$ gives a morphism
\[ \ome(\dpu) \mto{\phi} \ome \]
equivariant under the symplectic group Sp$(W)$. Since the image of 
$\phi$ is $k \oplus W \oplus \im \phi_2$, $T(\phi)$ corresponds to 
a null correlation bundle on $\proj{W}$ (which is Sp$(W)$-equivariant).


\begin{thebibliography}{XXXX}


\bibitem{BGG} I. Bernstein, I. Gelfand, S. Gelfand. {\em Algebraic bundles
over $\pn$ and problems of linear algebra.}
Funkts. Anal. Prilozh. {\bf 12} (1978); english translation in
Functional analysis and its applications {\bf 12} (1978), 212-214. 

\bibitem{Ei} D. Eisenbud. {\em Commutative algebra.} 
GTM 150, Springer-Verlag (1995).

\bibitem{EFS} D. Eisenbud, G. Fl\o ystad, F.-O. Schreyer. {\em
Sheaf Cohomology and Free Resolutions over Exterior Algebras.}
 {\tt http:// arXiv.org/abs/math.AG/0104203.}

\bibitem{Fl} G. Fl\o ystad. {\em Describing coherent sheaves on 
projective spaces via Koszul duality.} 
 {\tt http:// arXiv.org/abs/math.AG/0012263.}

\bibitem{Ja} J.C. Jantzen. {\em Representations of algebraic groups.}
Pure and Applied Mathematics {\bf 131} (1987), Academic Press Inc.

\end{thebibliography}
\end{document}